\documentclass{amsart}
\usepackage{amsmath,latexsym}
\usepackage[mathscr]{eucal}
\usepackage{amssymb}

\allowdisplaybreaks 
\theoremstyle{remark}

\def\cite#1{{\rm [#1]}}
\begin{document}

\title{The path to recent progress on small gaps between primes}
\author{D. A. Goldston, J. Pintz and C. Y. Y{\i}ld{\i}r{\i}m}

\thanks{The first author was supported by NSF grant DMS-0300563, the
NSF Focused Research Group grant 0244660, and the American Institute of
Mathematics; the second author by OTKA grants No. T38396, T43623, T49693 and
the Balaton program; the third author by T\"{U}B\.{I}TAK }

\date{\today}

\maketitle

\section*{1. Introduction }
In the articles {\it Primes in Tuples I \& II} ([13], [14])
we have presented the
proofs of some assertions about the existence of small gaps between
prime numbers which go beyond the hitherto established results. Our method
depends on tuple approximations. However, the approximations and the way
of applying the approximations has changed over time, and some comments in
this paper may provide insight as to the development of our work. 

First, here is a short narration of our results. Let
\begin{equation}
\theta(n) := \begin{cases}
\log n &\text{ if $n$ is prime}, \\
0 &\text{ otherwise},
\end{cases}
\label{eq: 1}
\end{equation}
and
\begin{equation}
\Theta(N;q,a) := \sum_{\substack{n\leq N \\ n\equiv a \, (\bmod
\, q)}}\theta(n).
\label{eq: 2}
\end{equation}
In this paper $N$ will always be a large integer,
$p$ will denote a prime number, and $p_{n}$ will denote the $n$-th prime.
The prime number theorem says that 
\begin{equation}
\lim_{x\to\infty}{|\{p: \; p \leq x\}|\over {x\over \log x}} 
=1 ,
\label{eq: 3}
\end{equation}
and this can also be expressed as
\begin{equation}
\sum_{n \leq x}\theta(n) \sim x  \qquad {\rm as} \,\, x\to\infty.
\label{eq: 4}
\end{equation}
It follows trivially from the prime
number theorem that
\begin{equation}
\liminf_{n\to \infty}{p_{n+1}-p_{n}\over \log p_{n}} \leq 1 .
\label{eq: 5}
\end{equation}
By combining former methods with
a construction of certain (rather sparsely distributed) intervals
which contain more primes than the expected number by a factor of
$e^{\gamma}$, Maier [25] had reached the best known result in this direction
that
\begin{equation}
\liminf_{n\to \infty}{p_{n+1}-p_{n}\over \log p_{n}} \leq 0.24846... \;\; .
\label{eq: 6}
\end{equation}

It is natural to expect that modulo $q$ the primes would be almost equally
distributed in the reduced residue classes. The deepest knowledge
on primes which plays a role in our method concerns a measure of
the distribution of primes in reduced residue classes referred to as the 
level of distribution of primes in arithmetic progressions.
We say that the
primes
have {\it level of distribution} $\alpha$ if
\begin{equation}
\sum_{q\leq Q}\max_{\substack{a \\ (a,q)=1}}\left|
\Theta(N;q,a) - {N\over \phi(q)}\right| \ll {N\over (\log N)^{A}}
\label{eq: 7}
\end{equation}
holds for any $A > 0$ and any arbitrarily small fixed $\epsilon > 0$ with
\begin{equation}
Q = N^{\alpha - \epsilon} .
\label{eq: 8}
\end{equation}
The {\it
Bombieri-Vinogradov theorem} provides the level ${1\over 2}$,
while the {\it Elliott-Halberstam conjecture} asserts that the primes have
level of distribution $1$.

The Bombieri-Vinogradov theorem
allows taking $Q=N^{{1\over 2}}(\log N)^{-B(A)}$ in (7), by virtue of which
we have proved unconditionally in [13] that
for any fixed $r\geq 1$,
\begin{equation}
\liminf_{n\to \infty}{p_{n+r}-p_{n}\over \log p_{n}} \leq (\sqrt{r}-1)^2
\;\,  ;
\label{eq: 9}
\end{equation}
in particular,
\begin{equation}
\liminf_{n\to \infty}{p_{n+1}-p_{n}\over \log p_{n}} = 0 .
\label{eq: 10}
\end{equation}
In fact, assuming that the level of distribution of primes is $\alpha$, we 
obtain more generally than (9) that, for $r \geq 2$, 
\begin{equation}
\liminf_{n\to \infty}{p_{n+r}-p_{n}\over \log p_{n}} \leq (\sqrt{r}-
\sqrt{2\alpha})^2 .
\label{eq: 11}
\end{equation}
Furthermore, assuming that $\alpha > 
{1\over 2}$, there exists an explicitly calculable constant $C(\alpha)$
such that for $k \geq C(\alpha)$ any sequence of $k$-tuples 
\begin{equation}
\{(n+h_{1}, n+h_{2}, \ldots, n+h_{k})\}_{n=1}^{\infty} ,
\label{eq: 12}
\end{equation}
with the set of distinct integers 
$\mathcal H = \{h_{1}, h_{2}, \ldots, h_{k}\}$ 
{\it admissible} in the sense
that $\displaystyle \prod_{i=1}^{k}(n+h_{i})$ has no fixed prime factor
for every $n$, contains at least two primes infinitely often. For instance
if $\alpha \geq 0.971$, then this holds for $k \geq 6$, giving
\begin{equation}
\liminf_{n\to \infty} (p_{n+1}-p_{n}) \leq 16,
\label{eq: 13}
\end{equation}
in view of the shortest admissible 
$6$-tuple $(n, n+4, n+6, n+10, n+12, n+16)$.

We note that the gaps obeying Eq.s (9)-(11) constitute a positive proportion
of all gaps of the corresponding kind. By incorporating Maier's method
into ours we improved (9) to
\begin{equation}
\liminf_{n\to \infty}{p_{n+r}-p_{n}\over \log p_{n}} \leq 
e^{-\gamma}(\sqrt{r}-1)^2 ,
\label{eq: 14}
\end{equation}
but for these gaps we don't have a proof of there being a
positive proportion of all gaps of this kind. (These results will appear in
forthcoming articles).

In [14] the result (10) was considerably improved to
\begin{equation}
\liminf_{n\to \infty}{p_{n+1}-p_{n}\over (\log p_{n})^{{1\over 2}}
(\log\log p_{n})^2 } < \infty.
\label{eq: 15}
\end{equation}
In fact, the methods of [14] lead to a much more general result:
When $\mathcal A \subseteq \Bbb N$ is a sequence satisfying 
$\mathcal A (N) := |\{n; n\leq N, n\in \mathcal A \}| > C(\log
N)^{1/2}(\log\log
N)^2 $ for all sufficiently large $N$, infinitely many of the differences of
two elements of $\mathcal A$
can be expressed as the difference of two primes.

\section*{2. Former approximations by truncated divisor sums}

The von Mangoldt function 
\begin{equation}
\Lambda(n) := \begin{cases}
\log p &\text{ if \, $n=p^{m}$, $m\in \Bbb Z^{+}$}, \\
0 &\text{ otherwise},
\end{cases}
\label{eq: 16}
\end{equation}
can be expressed as
\begin{equation}
\Lambda(n) = \sum_{d \mid n }\mu(d)\log({R\over d}) \qquad {\rm for} \,\, n
> 1 .
\label{eq: 17}
\end{equation}
Since the proper prime powers contribute negligibly, the prime number
theorem (4) can be rewritten as
\begin{equation}
\psi(x) := \sum_{n \leq x}\Lambda(n) \sim x  \qquad {\rm as} \,\,
x\to\infty.
\label{eq: 18}
\end{equation}
It is natural to expect that the truncated sum
\begin{equation}
\Lambda_{R}(n) := \sum_{\substack{d \mid n \\ d \leq R}}\mu(d)\log({R\over
d}) \qquad {\rm for} \,\, n\geq 1.
\label{eq: 19}
\end{equation}
mimics the behaviour of $\Lambda (n)$ on some averages. 

The beginning of our line of research is Goldston's [6] alternative 
rendering of the proof of
Bombieri and Davenport's theorem on small gaps between primes. Goldston
replaced the application of the circle method in the original proof by the
use of the truncated divisor sum (19). The use of functions like
$\Lambda_{R}(n)$ 
goes back to Selberg's work [27] on the zeros of the Riemann zeta-function 
$\zeta(s)$. The most beneficial feature of the truncated divisor sums is that 
they can be used in place of $\Lambda(n)$ on some occasions when it is not 
known how to work with $\Lambda(n)$ itself. The principal such situation
arises in counting the primes in tuples. Let
\begin{equation}
\mathcal H = \{h_{1}, h_{2}, \ldots, h_{k}\} \;\; {\rm with} \;\; 1\leq
h_{1}, \ldots, h_{k} \leq h \;\; {\rm distinct \,\, integers} 
\label{eq: 20}
\end{equation}
(the restriction of $h_{i}$ to positive integers is inessential; the whole
set $\mathcal H$ can be shifted by a fixed integer with no effect on our
procedure), and for a prime $p$ denote by
$\nu_{p}(\mathcal H)$ the number of distinct residue classes
modulo $p$ occupied by the elements of $\mathcal H$. The singular series
associated with the $k$-tuple $\mathcal H$ is defined as
\begin{equation}
\mathfrak S (\mathcal H)
: = \prod_{p}(1-{1\over p})^{-k}(1-{\nu_{p}(\mathcal H)\over p}).
\label{eq: 21}
\end{equation}
Since $\nu_{p}(\mathcal H)=k$ for $p>h$, the product is convergent.
The admissibility of $\mathcal H$ is equivalent to $\mathfrak S
(\mathcal H) \neq 0$,
and to $\nu_{p}(\mathcal H) \neq p$ for all primes. Hardy and Littlewood
[23] conjectured that
\begin{equation}
\sum_{n \leq N}\Lambda(n;\mathcal H) :=
\sum_{n\leq N}\Lambda(n+h_{1}) \cdots \Lambda(n+h_{k}) = N(\mathfrak S
(\mathcal H) + o(1)), \quad {\rm as} \,\, N\to\infty .
\label{eq: 22}
\end{equation}
The prime number theorem is the $k=1$ case, and for $k\geq 2$ the conjecture
remains unproved.
(This conjecture is trivially true if $\mathcal H$ is inadmissible). 

A simplified version of Goldston's argument in [6] was given in [17] as
follows. To obtain information on small gaps between primes, let 
\begin{equation}
\psi(n,h) := \psi(n+h)-\psi(n,h) = \sum_{n<m\leq n+h}\Lambda(m), \quad
\psi_{R}(n,h) := \sum_{n<m\leq n+h}\Lambda_{R}(m) ,
\label{eq: 23}
\end{equation}
and consider the inequality
\begin{equation}
\sum_{N < n\leq 2N}(\psi(n,h)-\psi_{R}(n,h))^{2} \geq 0.
\label{eq: 24}
\end{equation}
The strength of this inequality depends on how well $\Lambda_{R}(n)$
approximates $\Lambda(n)$. On multiplying out the terms and using
from [6] the formulas
\begin{align}
& \sum_{n\leq N}\Lambda_{R}(n)\Lambda_{R}(n+k) \sim \mathfrak S(\{0,k\})N,
\quad
\sum_{n\leq N}\Lambda(n)\Lambda_{R}(n+k) \sim \mathfrak S(\{0,k\})N \quad
(k\neq 0) \label{eq: 25}  \\
& \sum_{n\leq N}\Lambda_{R}(n)^{2} \sim N\log R, \quad 
\sum_{n\leq N}\Lambda(n)\Lambda_{R}(n) \sim N\log R ,
\label{eq: 26}
\end{align}
valid for $|k| \leq R \leq N^{{1\over 2}}(\log N)^{-A}$,
gives, taking $h=\lambda\log N$ with $\lambda \ll 1$,
\begin{equation}
\sum_{N < n\leq 2N}(\psi(n+h)-\psi(n))^{2} \geq (hN\log R + Nh^{2})(1-o(1)) 
\geq ({\lambda\over 2}+ \lambda^2 - \epsilon)N(\log N)^2 
\label{eq: 27}
\end{equation}
(in obtaining this one needs the two-tuple case of Gallagher's singular
series average given in (46) below, which can be traced back to Hardy and
Littlewood's and Bombieri and Davenport's work).
If the interval $(n, n+h]$ never contains more than one prime, then the
left-hand side of (27) is at most
\begin{equation}
\log N \sum_{N < n\leq 2N}(\psi(n+h)-\psi(n)) 
\sim \lambda N (\log N)^2 ,
\label{eq: 28}
\end{equation}
which contradicts (27) if $\lambda > {1\over 2}$, and thus one obtains
\begin{equation}
\liminf_{n\to \infty}{p_{n+1}-p_{n}\over \log p_{n}} \leq {1\over 2} .
\label{eq: 29}
\end{equation}
Later on Goldston et al. in [3], [4], [7], [15], [16], [18] applied this 
lower-bound method to various problems concerning the distribution
of primes and in [8] to the pair correlation of zeros of the Riemann
zeta-function. In most of these works the more delicate divisor sum
\begin{equation}
\lambda_{R}(n) := \sum_{r\leq R}{\mu^{2}(r)\over \phi(r)}\sum_{d\mid
(r,n)}d\mu(d)
\label{eq: 30}
\end{equation}
was employed especially because it led to better conditional results 
which depend on the Generalized Riemann Hypothesis.

The left-hand side of (27) is the second moment for primes in short
intervals. Gallagher [5] showed that
the Hardy-Littlewood conjecture (22) implies that the
moments for primes in intervals of length $h \sim \lambda\log N$ are the
moments of a Poisson distribution with mean $\lambda$. In particular, 
it is expected that 
\begin{equation}
\sum_{n\leq N}(\psi(n+h)-\psi(n))^2 \sim (\lambda + \lambda^{2})
 N (\log N)^2 
\label{eq: 31}
\end{equation}
which in view of (28) implies (10) but 
is probably very hard to prove. It is known from the 
work of Goldston and Montgomery [12] that assuming the Riemann Hypothesis,
an extension of (31) for $1 \leq h \leq N^{1-\epsilon}$
is equivalent to a form of the pair correlation conjecture for the
zeros of the Riemann zeta-function. We thus see that the factor 
${1\over 2}$ in (27) is what is lost from the truncation level $R$, and an
obvious strategy is to try to improve on the range of $R$ where
(25)-(26) are valid. In fact, the asymptotics in (26) are
known to hold for $R\leq N$ (the first relation in (26) is a special case
of a result of Graham [21]).
It is easy to see that the second relation in
(25) will hold with $R=N^{\alpha - \epsilon}$, where $\alpha$ is the level
of distribution of primes in arithmetic progressions. For the first relation
in (25) however, one can prove the the formula is valid for $R=N^{1/2 +
\eta}$ for a small $\eta > 0$, but unless one also assumes a somewhat
unnatural level of distribution conjecture for $\Lambda_{R}$, one can go no
further. Thus increasing the range of $R$ in (25) is not currently possible.

However, there is another possible approach motivated by Gallagher's
work [5]. In 1999 the first and third authors discovered how to calculate
some of the higher moments of the short divisor sums (19) and (30). At first
this was
achieved through straightforward summation and only the triple correlations
of $\Lambda_{R}(n)$ were worked out in [17]. In applying these formulas, the
idea of finding approximate moments with some expressions corresponding to
(24) was eventually replaced with
\begin{equation}
\sum_{N < n\leq 2N}(\psi(n,h)-\rho\log N)(\psi_{R}(n,h)-C)^{2}
\label{eq: 32}
\end{equation}
which if positive for some $\rho > 1$ implies that for some $n$ we have
$\psi(n,h) \geq 2\log N $. Here $C$ is available to optimize the
argument. Thus the
problem was switched from trying to find a good fit for $\psi(n,h)$ with a
short divisor sum approximation to the easier problem of trying to maximize
a given quadratic form, or more generally a mollification problem. With
just third correlations this resulted in (29), thus giving no improvement
over Bombieri and Davenport's result. Nevertheless the new method was not
totally fruitless since it gave
\begin{equation}
\liminf_{n\to \infty}{p_{n+r}-p_{n}\over \log p_{n}} \leq 
r - {\sqrt{r}\over 2} ,
\label{eq: 33}
\end{equation}
whereas the argument leading to (29) gives $r-{1\over 2}$. Independently of
us, Sivak [29] incorporated Maier's method into
[17] and improved upon (33) by the factor $e^{-\gamma}$ (cf. (6) and
(14) ). 

Following [17],
with considerable help from other mathematicians, in [20] the
$k$-level correlations of $\Lambda_{R}(n)$ were calculated. This leap
was achieved through replacing straightforward summation with complex
integration upon the use of Perron type formulae. Thus it became feasible to
approximate $\Lambda(n,\mathcal H)$ which was defined in (22) by 
\begin{equation}
\Lambda_{R}(n;\mathcal H) := \Lambda_{R}(n +h_{1})\Lambda_{R}(n
+h_{2}) \cdots \Lambda_{R}(n +h_{k}) .
\label{eq: 34}
\end{equation}
Writing
\begin{equation}
\Lambda_{R}(n; \mathbf H) := (\log R)^{k-|\mathcal H|}\Lambda_{R}(n;\mathcal
H), \quad
\psi_{R}^{(k)}(n,h) := \sum_{1\leq h_{1},\ldots, h_{k}\leq h}
\Lambda_{R}(n; \mathbf H) ,
\label{eq: 35}
\end{equation}
where the distinct components of the $k$-dimensional 
vector $\mathbf H$ are the elements of the set $\mathcal H$,
$\psi_{R}^{(j)}(n,h)$ provided the approximation to $\psi(n,h)^j$, and 
the expression 
\begin{equation}
\sum_{N < n\leq 2N}(\psi(n,h)-\rho\log N)(\sum_{j=0}^{k}
a_{j}\psi_{R}^{(j)}(n,h)(\log R)^{k-j})^2
\label{eq: 36}
\end{equation}
could be evaluated. Here the $a_{j}$ are constants available to optimize the
argument. The optimization turned out to be a rather complicated problem
which will not be discussed here, but the solution was recently completed in
[19] with the result that for any fixed $\lambda >
(\sqrt{r}-\sqrt{{\alpha\over 2}})^2 $ and $N$ sufficiently large,  
\begin{equation}
\sum_{\substack{n\leq N \\ p_{n+r}-p_{n} \leq \lambda\log p_{n}}}1 \gg_{r}
\sum_{\substack{p\leq N \\ p :\, {\rm prime}}}1 .
\label{eq: 37}
\end{equation}
In particular, unconditionally,
for any fixed $\eta >0$ and for all sufficiently large $N >
N_{0}(\eta)$, a positive proportion of gaps $p_{n+1}-p_{n}$ with $p_{n}\leq
N$ are smaller than $({1\over 4}+\eta)\log N$.
This is numerically a little short of Maier's result (6), but (6) was shown
to hold for a sparse sequence of gaps. The work [19] also turned out to be
instrumental in Green and Tao's [22] proof that the primes contain
arbitrarily long arithmetic progressions.

The efforts made in 2003 using divisor sums which are more complicated than
$\Lambda_{R}(n)$ and $\lambda_{R}(n)$ gave rise to more difficult
calculations and didn't meet with success.
During this work Granville and Soundararajan provided us with the idea that
the method should be applied directly to individual tuples rather than sums
over tuples which constitute approximations of moments. They replaced
the earlier expressions with
\begin{equation}
\sum_{N < n\leq 2N}(\sum_{h_{i}\in \mathcal H}\Lambda(n+h_{i}) - r\log
3N)(\tilde{\Lambda}_{R}(n; \mathcal H))^2  ,
\label{eq: 38}
\end{equation}
where $\tilde{\Lambda}_{R}(n; \mathcal H)$ is a short divisor sum which
should be large when $\mathcal H$ is a prime tuple. This is the type of
expression which is used in the proof of the result described in connection
with (12)--(13) above. However, for obtaining the results (9)--(11).
we need arguments based on using (32) and (36).

\section*{3. Detecting prime tuples}

We call the tuple (12) a {\it prime tuple} when all of its components
are prime numbers. Obviously this is equivalent to requiring that 
\begin{equation}
P_{\mathcal H}(n) := (n+h_{1})(n+h_{2})\cdots (n+h_{k})
\label{eq: 39}
\end{equation}
is a product of $k$ primes. As the generalized von Mangoldt function
\begin{equation}
\Lambda_{k}(n) := \sum_{d \mid n}\mu(d)(\log {n\over d})^k
\label{eq: 40}
\end{equation}
vanishes when $n$ has more than $k$ distinct prime factors, we may use
\begin{equation}
{1\over k!}\sum_{\substack{d\mid
P_{\mathcal H}(n) \\ d \leq R}}\mu(d)(\log {R\over d})^k
\label{eq: 41}
\end{equation}
for approximating prime tuples. (Here $1/k!$ is just a normalization factor.
That (41) will be also counting some tuples by including proper prime 
power factors doesn't pose a threat since in our applications their
contribution is negligible). But this idea by itself
brings restricted progress: now the right-hand side of (6) can be replaced
with $1-{\sqrt{3}\over 2}$.

The efficiency of the argument is greatly increased
if instead of trying to include tuples composed only of primes, one looks
for tuples with primes in many components. So in [13] we employ
\begin{equation}
\Lambda_{R}(n; \mathcal H, \ell) := {1\over (k+\ell)!}\sum_{\substack{d\mid
P_{\mathcal H}(n) \\ d \leq R}}\mu(d)(\log {R\over d})^{k+\ell} ,
\label{eq: 42}
\end{equation}
where $|\mathcal H| = k$ and $0 \leq \ell \leq k$,
and consider those $P_{\mathcal H}(n)$ which have at most $k+\ell$ 
distinct prime factors. In our applications the optimal
order of magnitude of the integer $\ell$ turns out to be about $\sqrt{k}$. 
To implement this new approximation in the skeleton of the argument, the
quantities
\begin{equation}
\sum_{n \leq N}\Lambda_{R}(n;\mathcal H_{1},\ell_{1})
\Lambda_{R}(n;\mathcal H_{2},\ell_{2}) ,
\label{eq: 43}
\end{equation}
and
\begin{equation}
\sum_{n \leq N}\Lambda_{R}(n;\mathcal H_{1},\ell_{1})
\Lambda_{R}(n;\mathcal H_{2},\ell_{2}) \theta(n+h_{0}) ,
\label{eq: 44}
\end{equation}
are calculated as $R,\, N \to \infty$.
The latter has three cases according as $h_{0} \not\in
\mathcal H_{1} \cup \mathcal H_{2}$, or $h_{0} \in \mathcal H_{1} \setminus
\mathcal H_{2}$, or $h_{0} \in \mathcal H_{1} \cap \mathcal H_{2}$.
Here $M= |\mathcal H_{1}| + |\mathcal H_{2}| + \ell_{1} + \ell_{2} $
is taken as a fixed integer which may be arbitrarily large.
The calculation of (43) is valid with
$R$ as large as $N^{{1\over 2}-\epsilon}$ and $h \leq R^{C}$ for any
constant $C>0$.
The calculation of (44) can be carried out for
$R$ as large as $N^{{\alpha\over 2}-\epsilon}$ and $h \leq R$.
It should be noted that in [19] in the same context the usage of
(34) which has $k$ truncations,
restricted the range of the divisors greatly, for 
then $R \leq N^{{1\over 4k}-\epsilon}$ was needed. Moreover the calculations
were more complicated compared to the present situation of dealing with only
one truncation.

Requiring the positivity of the quantity
\begin{equation}
\sum_{n=N+1}^{2N}(\sum_{1 \leq h_{0}\leq h}\theta(n+h_{0}) - r\log
3N)(\sum_{\substack{\mathcal H \subset \{1,2,\ldots , h\} \\ |\mathcal
H|=k}}\Lambda_{R}(n; \mathcal H, \ell))^2 ,  \qquad (h=\lambda\log 3N),
\label{eq: 45}
\end{equation}
which can be calculated easily from asymptotic formulas for 
(43) and (44), and
Gallagher's [5] result that with the notation of (20) for fixed $k$ 
\begin{equation}
\sum_{\mathcal H}\mathfrak S (\mathcal H)
\sim h^k \qquad {\rm as} \;\;\; h\to\infty ,
\label{eq: 46}
\end{equation}
yields the results (9)--(11). For the proof of the result mentioned in 
connection with (12), the positivity of (38) with $r=1$ and
$\Lambda_{R}(n; \mathcal H,\ell)$ for an $\mathcal H$ satisfying (20) 
in place of $\tilde{\Lambda}_{R}(n; \mathcal H)$ is used. For (13),
the positivity of 
an optimal linear combination of the quantities for (12) is pursued.

The proof of (15) in [14] also depends on the positivity of (45) for $r=1$ 
and $h = {C\log N\over k}$ modified with the extra restriction
\begin{equation}
(P_{\mathcal H}(n), \displaystyle \prod_{p\leq \sqrt{\log N}}p)=1 
\label{eq: 47}
\end{equation}
on the tuples to be summed over,
but involves some essential differences from the procedure 
described above. Now the size of $k$ is taken as large as
$c{\sqrt{\log N}\over (\log\log N)^2 }$ (where $c$ is a sufficiently small
explicitly calculable absolute constant). This necessitates a
much more refined treatment of the error terms arising in the argument,
and in due course the restriction (47) is brought in to avoid the
complications arising from the possibly irregular behaviour of
$\nu_{p}(\mathcal H)$ for small $p$. In the new argument
a modified version of the Bombieri-Vinogradov theorem is needed. 
Roughly speaking, in the version developed for this purpose, compared to (7) the 
range of the moduli $q$ is curtailed a little bit in return for a little
stronger upper-bound. Moreover, instead of Gallagher's result (46) which
was for fixed $k$ (though the result may hold for $k$ growing as some 
function of $h$, we do not know exactly how large this function can be
in addition to dealing with the problem of non-uniformity in $k$),
the weaker property that
$\sum_{\mathcal H}\mathfrak S (\mathcal H)/h^k$ is non-decreasing 
(apart from a factor of $1+o(1)$) as a function of $k$ 
is proved and employed. The whole argument is
designed to give the more general result which was mentioned after (15).

\section*{4. Small gaps between almost primes}

In the context of our work by {\it almost prime} we mean an {\it
$E_2$-number}, i.e. a natural number
which is a product of two distinct primes. 
We have been able to apply our methods to finding small gaps between almost
primes in collaboration with S. W. Graham. For this purpose a
Bombieri-Vinogradov type theorem for $\Lambda \ast \Lambda$ is needed, and
the work of Motohashi [26] on obtaining such a result for the Dirichlet
convolution of two sequences is readily applicable (see also [1]). In [9]
alternative proofs of some results of [13] such as (10) and (13) are
given couched in the formalism of the Selberg sieve. Denoting
by $q_{n}$ the $n$-th $E_{2}$-number, in [9] and [10] it is shown that
there is a constant $C$ such that for any positive integer $r$,
\begin{equation}
\liminf_{n\to\infty}(q_{n+r}-q_{n}) \leq Cre^{r} ;
\label{eq: 48}
\end{equation}
in particular
\begin{equation}
\liminf_{n\to\infty}(q_{n+1}-q_{n}) \leq 6 .
\label{eq: 49}
\end{equation}
Furthermore in [11] proofs of a strong form of the Erd\"{o}s--Mirsky
conjecture and related assertions have been obtained. 

\section*{5. Further remarks on the origin of our method}

In 1950 Selberg was working on applications of his sieve method to the twin
prime and Goldbach problems and invented a weighted sieve method that gave
results which were later superseded by other methods and thereafter
largely neglected. Much later in 1991 Selberg published the details of this
work in Volume II of his Collected Works [28], describing it as ``by now of
historical interest only". In 1997 Heath-Brown [24] generalized Selberg's
argument from the twin prime problem to the problem of almost prime tuples.
Heath-Brown let
\begin{equation}
\Pi = \prod_{i=1}^{k}(a_{i}n +b_{i})
\label{eq: 50}
\end{equation}
with certain natural conditions on the integers $a_{i}$ and $b_{i}$. Then
the argument of Selberg (for the case $k=2$) and Heath-Brown for the 
general case is to choose $\rho > 0$ and the numbers $\lambda_{d}$ of the
Selberg sieve so that, with $\tau$ the divisor function,
\begin{equation}
Q = \sum_{n\leq x}\{1-\rho\sum_{i=1}^{k}\tau(a_{i}n + b_{i})\}(\sum_{d \mid
\Pi}\lambda_{d})^2 > 0 .
\label{eq: 51}
\end{equation}
From this it follows that there is at least one value of $n$ for which
\begin{equation}
\sum_{i=1}^{k}\tau(a_{i}n + b_{i}) < {1\over \rho} .
\label{eq: 52}
\end{equation}
Selberg found in the case $k=2$ that $\rho = {1\over 14}$ is acceptable,
which shows that one of $n$ and $n+2$ has at most two, while the other has
at most three prime factors for infinitely many $n$. Remarkably, this is
exactly the same type of tuple argument of Granville and Soundararajan which
we have used, and the similarity doesn't end here. Multiplying out, we have
$Q = Q_{1} - \rho Q_{2}$ where
\begin{equation}
Q_{1} = \sum_{n\leq x}(\sum_{d \mid \Pi}\lambda_{d})^2 > 0 , \qquad
Q_{2} = \sum_{i=1}^{k} \sum_{n\leq x}
\tau(a_{i}n + b_{i})\}(\sum_{d \mid \Pi}\lambda_{d})^2 > 0 .
\label{eq: 53}
\end{equation}
The goal is now to pick $\lambda_{d}$ optimally. As usual, the $\lambda_{d}$
are first made $0$ for $d>R$. At this point it appears difficult to find the
exact solution to this problem. Further discussion of this may be found in 
[28] and [24]. Heath-Brown, desiring to keep $Q_{2}$ small, made the
choice
\begin{equation}
\lambda_{d} = \mu(d)({\log (R/d)\over \log R})^{k+1},
\label{eq: 54}
\end{equation}
and with this choice we see
\begin{equation}
Q_{1} = {((k+1)!)^{2}\over (\log R)^{2k+2}}
\sum_{n\leq x}(\Lambda_{R}(n; \mathcal H, 1))^2 .
\label{eq: 55}
\end{equation}

Hence Heath-Brown used the approximation for a $k$-tuple with at most $k+1$
distinct prime factors. This observation was the starting point for our work
with the approximation $\Lambda_{R}(n; \mathcal H, \ell)$. The evaluation
of $Q_{2}$ with its $\tau$ weights is much harder to evaluate than $Q_{1}$
and requires Kloosterman sum estimates. The weight $\Lambda$ in $Q_{2}$ in
place of $\tau$ requires essentially the same analysis as
$Q_{1}$ if we use the Bombieri-Vinogradov theorem.
Apparently these arguments were never viewed as directly applicable
to primes themselves, and this connection was missed until now.

\vspace*{.5cm}
\footnotesize
D. A. Goldston  \,\,\, (goldston@math.sjsu.edu)

Department of Mathematics

San Jose State University 

San Jose, CA 95192

 USA \\

J. Pintz \,\,\, (pintz@renyi.hu)

R\'enyi Mathematical Institute of the Hungarian Academy of Sciences 

H-1364 Budapest 

P.O.B. 127 

Hungary \\

C. Y. Y{\i}ld{\i}r{\i}m \,\,\, (yalciny@boun.edu.tr)
\vspace*{-.1cm}
\begin{tabbing}
tw \= Department of Mathematicsmath \= \& math \=
\c Cengelk\"oy, Istanbul, P.K. 6, 81220 \=     \kill
\mbox{ } \>
Department of Mathematics  \> \mbox{ } \> Feza G\"ursey Enstit\"us\"u \\
\mbox{ } \> Bo\~{g}azi\c{c}i University \> \& \>
\c Cengelk\"oy, Istanbul, P.K. 6, 81220 \\
\mbox{ } \> Bebek, Istanbul 34342 \> \mbox{ } \> Turkey \\ 
\mbox{ } \> Turkey \> \mbox{ } \> \mbox{ }
\end{tabbing}

\end{document}